\documentclass[12pt]{amsart}
\usepackage{graphicx,amssymb}
\usepackage[all]{xy}
\usepackage{fullpage}
\newtheorem{thm}{Theorem}[section]
\newtheorem{cor}[thm]{Corollary}
\newtheorem{lem}[thm]{Lemma}
\newtheorem{prop}[thm]{Proposition}
\theoremstyle{definition}

\theoremstyle{remark}
\newtheorem{rem}[thm]{Remark}
\numberwithin{equation}{section}

\newcommand{\id}{\text{ id }}
 \DeclareMathOperator{\End}{End}
\DeclareMathOperator{\Gal}{Gal} 
\DeclareMathOperator{\Irr}{Irr}
\DeclareMathOperator{\Trace}{Trace}
\DeclareMathOperator{\Image}{Image}
\DeclareMathOperator{\Res}{Res}
\begin{document}

\title[Jacobians and rational idempotents]{Jacobians with group actions and rational idempotents}%
\author{Angel Carocca}
\address{Facultad de Matem\'aticas, Pontificia Universidad Ca\-t\'o\-li\-ca
de Chi\-le, Casilla 306--22, Santiago, Chile}
\email{acarocca@mat.puc.cl}
\author{Rub\'\i\ E. Rodr\'\i guez}
\address{Facultad de Matem\'aticas,  Pontificia Universidad Cat\'olica
de Chi\-le,  Casilla 306--22,  Santiago,  Chile}
\email{rubi@mat.puc.cl}
\thanks{This research was partially supported by
FONDECYT Grants \# 1030595 and \# 1011039, and the Presidential
Science Chair on Geometry Award.\\
The second author thanks the Institute for Mathematical Sciences,
SUNY at Stony Brook, for its generous hospitality while this work
was completed}

\maketitle
\thispagestyle{empty} \def\IMSmarkvadjust{0 pt}
\def\IMSmarkhadjust{0 pt}
\def\IMSmarkhpadding{0 pt}
\def\IMSpubltext{Published in modified form:}
\def\SBIMSMark#1#2#3{
 \font\SBF=cmss10 at 10 true pt
 \font\SBI=cmssi10 at 10 true pt
 \setbox0=\hbox{\SBF \hbox to \IMSmarkhpadding{\relax}
                Stony Brook IMS Preprint \##1}
 \setbox2=\hbox to \wd0{\hfil \SBI #2}
 \setbox4=\hbox to \wd0{\hfil \SBI #3}
 \setbox6=\hbox to \wd0{\hss
             \vbox{\hsize=\wd0 \parskip=0pt \baselineskip=10 true pt
                   \copy0 \break%
                   \copy2 \break% 
                   \copy4 \break}}
 \dimen0=\ht6   \advance\dimen0 by \vsize \advance\dimen0 by 8 true pt
                \advance\dimen0 by -\pagetotal
	        \advance\dimen0 by \IMSmarkvadjust
 \dimen2=\hsize \advance\dimen2 by .25 true in
	        \advance\dimen2 by \IMSmarkhadjust

%
%   Check for publication info
%
%  \newread\jref
  \openin2=publishd.tex
  \ifeof2\setbox0=\hbox to 0pt{}
  \else 
     \setbox0=\hbox to 3.1 true in{
                \vbox to \ht6{\hsize=3 true in \parskip=0pt  \noindent  
                {\SBI \IMSpubltext}\hfil\break
                \input publishd.tex 
                \vfill}}
  \fi
  \closein2
  \ht0=0pt \dp0=0pt
 \ht6=0pt \dp6=0pt
 \setbox8=\vbox to \dimen0{\vfill \hbox to \dimen2{\copy0 \hss \copy6}}
 \ht8=0pt \dp8=0pt \wd8=0pt
 \copy8
 \message{*** Stony Brook IMS Preprint #1, #2. #3 ***}
}

\def\IMSmarkvadjust{-30pt}
\SBIMSMark{2003/01}{May 2003}{}

% ----------------------------------------------------------------
\tableofcontents

\section{Introduction}
\label{sec:intro} It is known that every finite group acts on some
curve, hence on some Jacobian, and therefore on some principally
polarized abelian variety (\textit{p.p.a.v.}). On the other hand,
every action of a finite group on an abelian variety induces an
isogeny decomposition into factors related to the rational
irreducible representations of the given group.

However, the existence of this decomposition is obtained from the
theory of group representations, and there has been no general
geometric interpretation of the factors.

The object of this paper is to prove some general results about
rational idempotents for a finite group $G$ and deduce from them
geometric information about the components that appear in the
decomposition of the Jacobian variety $JW$ of a curve $W$ with
$G-$action.

The simplest case of such a decomposition is when the group $G$ is
the group with two elements $G \cong \mathbb{Z}/2\mathbb{Z} =
\langle j : j^2 \rangle$ acting on a curve $W$. It is clear that
then the Jacobian variety $JW$ of $W$ also has an involution $j$
acting on it and, furthermore, that $JW$ is isogenous to the
product $B_1 \times B_2$, where $B_1 = \Image (1+j)$ is the
$G-$invariant part of $JW$ and $B_2 = \Image(1-i)$ is the
anti-invariant part. This decomposition was already observed by
Wirtinger in \cite{w} and used by Schottky and Jung in \cite{s-j}.
Observe that $B_1$ is isogenous to the Jacobian $JW_G$ of the
quotient curve $W_G = W/G$ and that $B_2$ was later called by
Mumford \cite{mum} the Prym variety $P(W/W_G)$ of the given cover.

The decomposition of Jacobians (or, more generally, abelian
varieties) with group actions has been studied in different
settings, with applications to theta functions, to the theory of
integrable systems and to the moduli spaces of principal bundles
of curves. For other special groups, there are the dihedral group
$D_p$ with $p$ prime in \cite{ries}, the symmetric group $S_3$ in
\cite{rr1}, all the subgroups of $S_4$ in \cite{rr2}, the
alternating group $A_5$ in \cite{sa1}, some Weyl groups in
\cite{kanev}, $S_4$ and $WD_4$ in \cite{d-m}, the dihedral groups
$D_n$ in \cite{crr} and $S_5$ in \cite{sa2} and \cite{lr2}.

For all these cases except the last two, it was shown that the
factors appearing in the decomposition were all of the form of
either Jacobians or Prym varieties of intermediate covers. In the
last two some of the factors were not of this type, but they could
still be identified as (connected components of the origin of)
intersections of Prym varieties of intermediate covers.

M\'erindol in \cite{m} (for abelian varieties) and Donagi in
\cite{d} (for Jacobians) gave a decomposition in the case when all
the irreducible rational representations of the group are
absolutely irreducible: i.e., when they stay irreducible when
considered as complex representations. In \cite{L-R}, Lange and
Recillas have recently given the general decomposition for any
abelian variety admitting a group action.

In the case studied in \cite{m} and \cite{d}, there are formulae
giving rational projectors whose images are the required factors
in the decomposition. These formulae are constructed based on the
fact that the condition imposed on the group $G$ implies the
existence of an essentially unique $G-$invariant inner product on
each of the vector spaces corresponding to the irreducible
representations, and therefore cannot be generalized.

In this work we find the corresponding rational projectors for the
general case, as well as describe rational projectors invariant
under any given subgroup. From them we deduce the decomposition of
any Prym or Jacobian variety of an intermediate cover, in the case
of a Jacobian with group action.

These explicit constructions allow geometric descriptions of the
factors appearing in the decomposition of a Jacobian with group
action. For instance, we give a necessary and sufficient condition
for a Prym variety of an intermediate cover to be such a factor.

In Section \ref{sec:not} we recall some basic results on
irreducible representations and idempotents in order to fix the
notation.

We then describe two methods for the construction of rational
idempotents of a finite group $G$: one for primitive idempotents
in Section \ref{sect:minimal} and another for idempotents that are
bilaterally invariant under any given subgroup of $G$ in Section
\ref{sec: inv}.

In Section \ref{sec:abvar}  we apply our constructions to the
decomposition of abelian varieties and Jacobians with group
action.

The paper concludes with two appendices, where we include two
examples; both groups have rational irreducible representations
which are \textit{not} absolutely irreducible and hence could not
have been studied without this general theory. The first one
illustrates our constructions of rational idempotents and the
second one exhibits a factor in the decomposition of the
corresponding Jacobians which is not of the previous known types;
that it, it is neither a Jacobian nor a Prym nor an intersection
of Prym varieties of intermediate covers.

\section{Preliminaries}
\label{sec:not}

Let $G$ be a finite group. With relation to representations, we
follow the notations in \cite{serre} and \cite{curtis-reiner}. The
known results quoted in this section may be found there as well.

We denote by $\Irr_F (G)$ the set of irreducible representations
of $G$ over the field $F$.

If $V \in \Irr_{\mathbb{C}}(G)$,  we let $L$ denote the field of
definition of $V$ and let $K$ denote the field obtained by
adjoining to the rational numbers $\mathbb{Q}$ the values of the
character $\chi_V$; then $K \subseteq L$ and $m =
m_{\mathbb{Q}}(V) = [L:K]$ is the Schur index of $V$.

If $H$ is any subgroup of $G$, $\rho_H$ will denote the
representation of $G$ induced by the trivial representation of
$H$. If we denote by $\langle U , V \rangle$ the usual inner
product of the characters of the representations $U$ and $V$, and
by $V^H$ the subspace of $V$ fixed under $H$, then it follows
immediately from the Frobenius reciprocity Theorem that
$$ \langle \rho_H , V \rangle = \dim_{\mathbb{C}} V^H \ .
$$

It is known that if $\Gal (L/\mathbb{Q})$, $\Gal (L/K)$ and
$\Gal(K/\mathbb{Q})$ denote the respective Galois groups, then
each representation $V^{\sigma}$  conjugate to $V$ by an element
$\sigma$ in $\Gal(L/\mathbb{Q})$ is also defined over $L$ and both
$V$ and $V^{\sigma}$ share the same field $K$. Furthermore,
$V^{\sigma}$ is equivalent to $V$ if and only if $\sigma$ is in
$\Gal (L/K)$.

We denote by $U$ the representation
$$ U = \bigoplus_{\tau \in \Gal (L/K)} V^{\tau} \simeq m \, V
$$

\noindent and remark that it is actually defined over $K$. The
Schur index $m$ of $V$ is minimal with respect to this property.

It is also known that the irreducible rational representation
$\mathcal{W}$ of $G$ associated to $V$ is characterized by its
decomposition into irreducible components over the corresponding
fields as follows.

\begin{equation}
\label{eq:WK}
 K \otimes_{\mathbb{Q}} \mathcal{W} \simeq \bigoplus_{\varphi
\in \Gal(K/\mathbb{Q})} U^{\varphi} \ , \text{ and }
\end{equation}

\begin{equation}
\label{eq:WL}
 L \otimes_{\mathbb{Q}} \mathcal{W} \simeq \bigoplus_{\sigma
\in \Gal(L/\mathbb{Q})} V^{\sigma} \simeq \bigoplus_{\varphi \in
\Gal(K/\mathbb{Q})} (m \, V)^{\varphi}
\end{equation}

The (unique) central idempotent $e_V$ of $L[G]$ that generates the
simple subalgebra associated to $V$ and the (unique) central
idempotent $e_{\mathcal{W}}$ of $\mathbb{Q}[G]$ that generates the
simple subalgebra associated to $\mathcal{W}$ are given as
follows.

$$ e_V = \frac{\dim V}{|G|} \sum_{g \in G} \chi_V (g^{-1}) \, g \
, \text{ and }
$$

\begin{equation}
\label{eq:ew} e_{\mathcal{W}} = \frac{\dim V}{|G|} \sum_{g \in G}
\Trace_{K/\mathbb{Q}}(\chi_V (g^{-1})) \, g \ .
\end{equation}

Note that $e_V$ is in fact an element of $K[G]$ and that it
coincides with the (unique) central idempotent $e_U$ of $K[G]$
that generates the simple subalgebra associated to $U = m \, V$;
furthermore, the following relations hold.

$$ e_V = e_{V^{\tau}} \ \text{ for all } \tau \in \Gal(L/K)
$$
and
$$ e_U e_{U^{\varphi}} = 0 \ \text{ for all } \varphi \in
\Gal(K/\mathbb{Q})\ .
$$

\medskip

It is also known that the simple algebras $L[G] \, e_V$,  $K[G] \,
e_V$ and $\mathbb{Q}[G] \, e_{\mathcal{W}}$ may be decomposed into
the direct sum of $n$, $\frac{n}{m}$ and $\frac{n}{m}$ minimal
left ideals, respectively. In other words, there exist primitive
idempotents $\{ \ell_1 , \ldots , \ell_n \}$, $\{k_1 \ldots ,
k_{\frac{n}{m}} \}$ and $\{ f_1 , \ldots , f_{\frac{n}{m}} \}$  in
the respective subalgebras, orthogonal among themselves, such that
the following relations hold.

\begin{equation}
\label{eq:ell}
 e_V = \ell_1 + \ldots + \ell_n \ \ \text{in } L[G] \, e_V
\end{equation}

$$ e_V = k_1 + \ldots + k_{\frac{n}{m}} \ \ \text{in } K[G] \, e_V
$$

$$ e_{\mathcal{W}} = f_1 + \ldots + f_{\frac{n}{m}} \ \ \text{in } \mathbb{Q}[G] \,
e_{\mathcal{W}}
$$

\medskip

These primitive idempotents are far from being unique, and except
for special cases, there are no general formulas for writing down
neither the $k_j$ nor the $f_j$.

A known case (\cite{m}) is the following. Suppose that the
rational irreducible representation $\mathcal{W}$ is
\textit{absolutely irreducible}; that is, it stays irreducible
when considered as a complex representation. Then $L = K =
\mathbb{Q}$, $m = 1$ and $n = \frac{n}{m} = \dim V = \dim
\mathcal{W}$.

In this case one can find the $\ell_j = k_j = f_j$'s as follows.
The hypothesis on $\mathcal{W}$ implies that there is a unique (up
to multiplication by a positive scalar) $G-$invariant inner
product on $\mathcal{W}$, denoted by $\langle \, , \, \rangle$.
Let $\{ w_1 , \ldots , w_n \}$ denote an orthogonal basis on
$\mathcal{W}$ and define
 $$ \ell_j =  \frac{n}{|G| \, ||w_j||^2} \sum_{g \in G} \langle w_j ,
 g(w_j) \rangle \, g \ .
 $$

To show that these $\ell_j$'s satisfy the required properties,
they are rewritten as follows

\begin{equation}
\label{eq:special} \ell_j =  \frac{\dim V}{|G|} \sum_{g \in G}
r_{jj} (g^{-1}) \, g \ , \text{ for } j \in \{ 1, \ldots , n \} \,
,
\end{equation}
where $\left( r_{ik} (g) \right)$ are the (rational) coefficients
of the matrix of the element $g$ with respect to the given basis.
Then (\ref{eq:ell}) is clear and that these $\ell_j$'s are
orthogonal idempotents is a consequence of the orthogonality
relations for characters (see \cite{serre}).

For the general case of any complex irreducible representation $V$
of $G$, of dimension $n$, one may \textit{define} $\ell_j$'s by
the right hand side of (\ref{eq:special}), where the $r_{ik}$'s
are now the coefficients of the given $L$-representation of $G$.
The same orthogonality relations show that these $\ell_j$'s are
primitive orthogonal idempotents in $L[G]$ whose sum is $e_V$.

In this paper we give a general construction for $k_j$'s and
$f_j$'s.

\section{The construction of primitive rational idempotents}
\label{sect:minimal}

Let $V$ be a complex irreducible representation of a finite group
$G$, with associated fields $L$ and $K$ as in Section
\ref{sec:not}.

\begin{prop}
\label{prop:stab}

Let $V$ be in $\Irr_{\mathbb{C}} G$.

Denote by $n = \dim V$, by $e = e_V$ the corresponding central
idempotent in $L[G]$, and by $\Gal (L/K) = \{ \tau_1 = 1, \ldots ,
\tau_m \} $, with $m$ the Schur index of $V$.

Consider any $n$ primitive orthogonal idempotents $\ell_1 , \ldots
, \ell_n $ in the simple algebra $L[G]\, e_V$ satisfying
(\ref{eq:ell}). For instance, the $\ell_j$ given by
(\ref{eq:special}).

Let $M_j = \sum_{h=1}^m L[G] \, \tau_h (\ell_j) $ be the sum of
the corresponding left ideals.

Then the stabilizer in $\Gal(L/K)$ of each left ideal $L[G] \,
\ell_j$ is trivial and the sum in $M_j$ is direct for each $j$ in
$\{ 1, \ldots, n \}$.

Furthermore, the left ideals $L[G] \, \tau_h (\ell_j)$ in each
$M_j$ are permuted among themselves by the action of $\Gal (L/K)$
and for each pair $j, k$ in $\{ 1, \ldots, n \}$, either $M_j =
M_k$ or $M_j \cap M_k = \{ 0 \}$.
\end{prop}

\begin{proof}
It is clear that $\tau (\ell_j)$ is a primitive idempotent in
$L[G]\, e_V$ for each $\tau$ in $\Gal(L/K)$.

Let $H$ denote the stabilizer of $\ell_j$ in $\Gal(L/K)$, and let
$\{ \tau_1 = 1, \tau_2 , ...., \tau_r \}$ be a right transversal
of $ \; H \; $ in $\Gal(L/K)$.

Consider
$$ M = L[G]\, \ell_j + \ldots + L[G] \, \tau_r(\ell_j) \ .
$$

Then $M$ is a left $K[G]$-module contained in $L[G]\, e_V $  with
$0 < \dim (M) \leq m \, \dim_L (V) = \dim_K (U)$.

But the minimality of $m$ then implies that $\dim(M) = m \, \dim
(V)$ and therefore $H = \{ 1\}$, $r = m$ and
$$ M = M_j = L[G]\, \ell_j \oplus \ldots \oplus L[G] \, \tau_m(\ell_j)$$

The other results are now clear, since $M_j \cap M_k$ is a left
$K[G]$-module contained in $L[G]\, e_V $  with dimension at most
$m \, \dim (V)$.
\end{proof}

\begin{rem}
At this point one could consider the elements of $K[G] \, e_V$
given by $\widetilde{k_j} = \sum_{h=1}^m \sigma_h (\ell_j) $ as
possible candidates to be the sought $k_j$. However, these
elements are not necessarily idempotents, since $\sigma_h
(\ell_j)$ and $\sigma_k (\ell_j)$ may be non orthogonal to each
other (e.g., see the example in Appendix \ref{sect:example}).
Furthermore, even if the $\widetilde{k_j}$ are generating
idempotents for $M_j$ for each $j$ in $\{ 1, \ldots, n \}$, there
would still remain the question of choosing $\frac{n}{m}$ among
them in such a way that they add up to $e_V$.
\end{rem}

Our next result gives an explicit way of finding $k_j$'s with all
the required properties.

\begin{thm}
\label{thm:kjs}

Let $V$ be in $\Irr_{\mathbb{C}} G$ and denote $n = \dim (V)$ and
$m = m_{\mathbb{Q}} (V)$.

Consider any $n$ primitive orthogonal idempotents $\ell_1 , \ldots
, \ell_n $ in the simple algebra $L[G]\, e_V$ satisfying
(\ref{eq:ell}).

Then there exist $\frac{n}{m}$ primitive orthogonal idempotents
$\{ u_{s}^1 \}_{s=1}^{\frac{n}{m}}$ in different left ideals $J_s
= L[G]\, \ell_{j_s}$ such that the following results hold.

\begin{enumerate}
 \item The element $u_s^h = \tau_h(u_s^1)$ is a
 primitive generating idempotent for the left
 ideal $\tau_h (J_s) = L[G] \, \tau_h(\ell_{j_s})$,
 for each $\tau_h$ in $\Gal(L/K)$ and for each $s$
 in $\{ 1, \ldots, \frac{n}{m} \}$.
 \item The $\{ u_s^h \}$ satisfy
  $$ u_s^h \, u_t^l =
\begin{cases}
    u_s^h, & \hbox{ if l=h and t=s;} \\
    0, & \hbox{otherwise.} \\
\end{cases}
  $$
\item The algebra $L[G] \, e_V$ decomposes as the direct
orthogonal sum of the left ideals generated by $\{ \sigma_h
(u_s^1) \}$ and furthermore we have
$$ e_V = \sum_{s = 1}^{\frac{n}{m}}
\sum_{h=1}^m \, u_s^h = \sum_{s = 1}^{\frac{n}{m}} \sum_{\sigma
\in \Gal(L/K)} \sigma(u_s^1) \ .
      $$
\end{enumerate}
\end{thm}

\begin{proof}
If we denote each minimal left ideal $L[G] \, \sigma_h (\ell_s)$
by $J_s^h$, it follows from Proposition \ref{prop:stab} that $M_1
= \sum_{h=1}^m J_1^h $ is a direct sum of the given left ideals.

If $m = n$, then $L[G] \, e_V = M_1$. If $m < n$, then there
exists $\ell_j$ with $j \neq 1$ such that $\ell_j \notin M_1$. By
renumbering the $\ell_j$ if needed, we may assume this is
$\ell_2$. Again from Proposition \ref{prop:stab} it follows that
$M_2$ is a direct sum of the corresponding ideals $J_2^h$ and,
furthermore, that the sum $M_1 + M_2$ of primitive left ideals $\{
J_s^h \}_{\substack{s=1,2 \\ h \in \{ 1, \ldots , m\} }}$ is also
direct.

If $M_1 + M_2$ is a proper subset of $L[G] \, e_V$, then we may
assume that for $\ell_3$ the corresponding $M_3$ has trivial
intersection with $M_1 + M_2$. Then, as before, the sum $M_1 + M_2
+ M_3$ of minimal left ideals $\{ J_s^h \}_{\substack{s=1,2,3 \\ h
\in \{ 1, \ldots , m\} }}$ is direct.

In any case, it is clear that the process terminates after
precisely $\frac{n}{m}$ steps, when we obtain the following direct
sum decomposition into minimal left ideals.

\begin{align}
\label{eq:directidealsum}
L[G] \, e_V & = J_1^1 \oplus J_1^2 \oplus \ldots \oplus J_1^m \\
\notag    & \oplus J_2^1 \oplus J_2^2 \oplus \ldots \oplus J_2^m \\
\notag    & \vdots \\
 \notag    & \oplus J_{\frac{n}{m}}^1 \oplus J_{\frac{n}{m}}^2 \oplus \ldots \oplus J_{\frac{n}{m}}^m \ .
\end{align}

Hence there exist unique primitive idempotents $u_s^h$ in $J_s^h$
such that
\begin{align}
\label{eq:idemp} e_V  & = u_1^1 + u_1^2 + \ldots + u_1^m \\
 \notag   & + u_2^1 + u_2^2 + \ldots + u_2^m \\
 \notag   & \vdots \\
 \notag   & + u_{\frac{n}{m}}^1 + u_{\frac{n}{m}}^2 + \ldots +
 \notag   u_{\frac{n}{m}}^m\ .
\end{align}

It follows that the $u_s^h$ are all orthogonal to each other, and
we claim that
 $$ u_s^h = \sigma_h (u_s^1)
 $$
for each $s$ and $h$.

Indeed, if we apply any $\tau_h \in \Gal(L/K)$ to the identity
(\ref{eq:idemp}), the left hand side is invariant and we obtain

\begin{align}
\notag e_V  & = \tau_h(u_1^1) + \tau_h(u_1^2) + \ldots + \tau_h(u_1^m) \\
 \notag   & + \tau_h(u_2^1) + \tau_h(u_2^2) + \ldots + \tau_h(u_2^m) \\
 \notag   & \vdots \\
 \notag   & + \tau_h(u_{\frac{n}{m}}^1) + \tau_h(u_{\frac{n}{m}}^2) + \ldots + \tau_h(u_{\frac{n}{m}}^m) \ .
\end{align}

\noindent where each $\tau_h(u_s^k)$ is a primitive idempotent in
some (unique) $J_s^{k_h} = \tau_h (\tau_k (J_s^1))$.

In particular, $\tau_h(u_s^1)$ is in $\tau_h(J_s^1) = J_s^h$ and
then the uniqueness of the decomposition (\ref{eq:idemp}) for the
given ideals (since the right hand side of
(\ref{eq:directidealsum}) is a direct sum) implies the desired
result.
\end{proof}

\begin{rem}
 \label{rem:const}
Starting from the $\ell_j$, the $u_s^h$ of Theorem \ref{thm:kjs}
may be found explicitly, possibly with some computer help. Our
example in Appendix \ref{sect:example} was built using the GAP
program.

The general algorithm is as follows, where we follow the notation
in the proof of the Theorem.

\begin{enumerate}
 \item Find a basis $B_s^h$ for each vector space $J_s^h$
 appearing in (\ref{eq:directidealsum}). It actually suffices to
 find a basis $B_s^1$ for each vector space $J_s^1$, since then a
 basis for each $J_s^h$ is given by $\tau_h(B_s^1)$.

\item Since  the sum on the right hand side of
(\ref{eq:directidealsum}) is direct, it follows that the union $B$
of all the bases $B_s^h$ is a basis for $L[G] \, e_V$.

Computing the coordinates of $e_V$ with respect to the basis $B$,
we find the desired primitive orthogonal idempotents $u_s^h$ in
each $J_s^h$.
\end{enumerate}

Alternatively, once the decomposition (\ref{eq:directidealsum})
has been obtained, one may use the method described in
\cite{fiedler} to actually find the $u_s^h$.
\end{rem}

\begin{cor}
\label{cor:ks}

Consider the $\{ u_s^1 \}_{s=1}^{\frac{n}{m}}$ constructed in
Theorem \ref{thm:kjs} and set
\begin{equation}
 \label{eq:ks}
k_s = \sum_{\tau \in \Gal(L/K)} \tau (u_s^1)
 \end{equation}
for each $s$ in $\{ 1, \ldots, \frac{n}{m} \}$.

Then $\{ k_s \}_{s=1}^{\frac{n}{m}}$ is a set of primitive
orthogonal idempotents in $K[G]$ such that

\begin{equation}
\label{eq:kjs}
 e_V = k_1 + \dots + k_{n/m} \ .
\end{equation}
\end{cor}

\begin{proof}
The results follow from the definition of the $k_s$ and the
properties of the $\tau (u_s^1)$ given by the Theorem.

In particular, the primitivity of each $k_s$ follows from the fact
that the dimension of the $L-$vector space $L[G] \, k_s = M_s$ is
equal to $m \, \dim_L (V) = \dim_K (U)$.
\end{proof}

\begin{cor}
\label{cor:fjs}

Given the idempotents $\ell_j$, construct the idempotents $k_s$ as
in Corollary \ref{cor:ks} and let $\mathcal{W}$ denote the
rational irreducible representation associated to $V$.

For each $s$ in $\{ 1, \ldots , \frac{n}{m} \}$ set

\begin{equation}
\label{eq:fs}
 f_s = \sum_{\varphi \in \Gal(K/\mathbb{Q})}\varphi(k_s)
     = \sum_{\sigma \in \Gal(L/\mathbb{Q})} \sigma(u_s^1) \ .
\end{equation}

Then the $f_s$ are primitive orthogonal idempotents in
$\mathbb{Q}[G]\, e_{\mathcal{W}}$ and they satisfy

\begin{equation}
\label{eq:efs} e_{\mathcal{W}} = f_1 + f_2 + \ldots +
f_{\frac{n}{m}} \ .
\end{equation}
\end{cor}

\begin{proof}
It is clear that each element $\varphi$ in $\Gal(K/\mathbb{Q})$
gives rise to primitive idempotents $\varphi(k_s)$, each
corresponding to the $K-$irreducible representation $U^{\varphi}$.

Since the $k_s$'s are orthogonal among themselves and their images
under different $\varphi$'s are associated to non conjugate
representations, we obtain the following relations for every
$\varphi_i , \varphi_j$ in $\Gal(K/\mathbb{Q})$ and for every $k_s
, k_t$.

$$ \varphi_i (k_s) \, \varphi_j (k_t) =
\begin{cases}
\varphi_i (k_s)& \text{if } j = i \text{ and } t = s;  \\
0 & \text{otherwise.}
\end{cases}
$$

Using these relations, it is easy to verify that the $f_s$'s are
indeed orthogonal idempotents. That their sum is $e_{\mathcal{W}}$
is an immediate consequence of (\ref{eq:kjs}) and (\ref{eq:ew}).

All that is left to show is that each $f_s$ is primitive, so
assume not; that is, assume there are rational orthogonal
idempotents $p_1$ and $p_2$ such that $f_s = p_1 + p_2$. But then
it follows from (\ref{eq:fs}) that $k_s = f_s \, e_U = p_1 \, e_U
+ p_2 \, e_U$, and this is a decomposition of the primitive
idempotent $k_s$ as the sum of two elements in $K[G] \, e_U$.

We will now prove that the $p_j \, e_U$ are orthogonal
idempotents, thus obtaining a contradiction. First, it is clear
that $p_i \, e_U \, p_j \, e_U = p_i \, p_j \, e_U$ and therefore
the $p_j \, e_U$ are orthogonal elements equal to their own
squares.

We just need to verify that they are not equal to zero. But if we
assume that $p_j \, e_U = 0$, we obtain

$$ p_j = p_j \, e_{\mathcal{W}} = \sum_{\varphi \in \Gal(K/\mathbb{Q})} \varphi(p_j \, e_U) = 0
$$

\noindent which is impossible.
\end{proof}

\section{Rational idempotents invariant under a subgroup}
\label{sec: inv}

In this section we find idempotents $f$ in $\mathbb{Q}[G]$ which
satisfy $h\, f = f \, h$ for all $h$ in a given subgroup $H$ of
$G$. This kind of idempotents will be needed in the study of the
decomposition of Jacobian varieties with $G-$action to be
developed in Section \ref{sec:Jac}.

\begin{rem}
\label{rem:turull}

Let $V$ denote a complex irreducible representation of a finite
group $G$ and let $K$ denote the field obtained by adjoining the
values of the character $\chi_V$ of $V$ to the rational numbers.

Let $H$ be a subgroup of $G$ and consider the central idempotent
$e_V = \frac{\dim V}{|G|} \sum_{g \in G} \chi_V(g^{-1}) \,  g$ in
$K[G]$ and the following elements
  \begin{align*}
  p_H & = \frac{1}{|H|} \sum_{h \in H} h \ \ , \text{ and }\\
   f_V & = p_H \, e_V \ \in K[G] \ .
  \end{align*}

Turull observes the following two facts in the proof of Theorem
1.4 of \cite{turull1}.

\begin{enumerate}
\item $p_H$ is an idempotent of $F[G]$ for any field $F$
containing $\mathbb{Q}$ and $F[G] \, p_H$ is an $F[G]-$module
affording the representation $\rho_H$ of $G$ induced by the
trivial representation of $H$, and
 \item if
 $$ \langle \rho_H \, , \, V  \rangle = 1
 $$
then $f_V$ is a primitive idempotent in $K[G]$ and $K[G] \, f_V $
is a $K[G]-$module affording $V$.
\end{enumerate}

We will see in Section \ref{sec:Jac} that (1) together with the
appropriate generalization of (2) given in Theorem \ref{thm:f}
turn out to be key points in understanding the isotypical
decomposition of a Jacobian variety with $G-$action.
\end{rem}

\begin{rem}
Note that Condition (2) above implies that the Schur index
$m_{\mathbb{Q}} (V)$ of $V$ is equal to one.

Since it is known that $m_{\mathbb{Q}} (V)$ divides $\langle
\rho_H \, , \, V \rangle$ for every subgroup $H$ of $G$ (see
\cite[Lemma 10.4]{isaacs}), a natural generalization of Condition
(2) would be to consider the case when
$$ \langle \rho_H \, , \, V  \rangle = m_{\mathbb{Q}} (V) \ ,
$$
which we do in Corollary \ref{cor:turullgen}.

The full general case
$$ \langle \rho_H \, , \, V  \rangle > 0
$$
is considered in Theorem \ref{thm:f}.
\end{rem}

The following elementary result will be used in the proof of
Theorem \ref{thm:f}. It is an immediate consequence of
(\ref{eq:WL}) and the Frobenius reciprocity theorem
(\cite{serre}), but the authors were unable to find it in the
literature.

\begin{lem}
\label{lem:genrho}

Denote $\Irr_{\mathbb{Q}} (G) = \{ \mathcal{W}_1 , \ldots ,
\mathcal{W}_r\}$ and for each $\mathcal{W}_j$ choose a
corresponding $V_j$ in $\Irr_{\mathbb{C}} (G)$. Also let $H$
denote a subgroup of $G$.

Then the rational isotypical decomposition of $\rho_H$ has the
following form

\begin{equation}
\label{eq:rho}
 \rho_H = a_1 \, \mathcal{W}_1 \oplus \ldots \oplus
a_r \, \mathcal{W}_r \ ,
\end{equation}

\bigskip

\noindent where
$$a_j = \frac{\langle \rho_H , V_j \rangle}{m_j} = \frac{\dim V_j^H}{m_j} \ .
$$
\end{lem}

\medskip

The following Theorem generalizes the second result in Remark
\ref{rem:turull}.

\begin{thm}
 \label{thm:f}
Let $\mathcal{W}$ be a rational irreducible representation of a
group $G$, and denote by $e_{\mathcal{W}}$ the associated central
idempotent in $\mathbb{Q}[G]$. We denote by $V$ one of the complex
irreducible representations of $G$ appearing in the isotypical
decomposition of $\mathbb{C} \otimes_{\mathbb{Q}} \mathcal{W}$.

For any subgroup $H$ of $G$, let
 $$ p_H = \frac{1}{|H|} \sum_{h \in H} h
 $$
be the central idempotent in $ \mathbb{Q}[H]$ corresponding to the
trivial representation of $H$.

Then
 $$ f_H =p_H \, e_{\mathcal{W}} = e_{\mathcal{W}} \, p_H
 $$
is an element of the simple algebra $\mathbb{Q}[G] \,
e_{\mathcal{W}}$ satisfying the following conditions.

\begin{enumerate}
\item $f_H^2 = f_H$,
 \item $ h \, f_H = f_H = f_H \, h $ for every $h$ in $H$, and
 \item $f_H = 0$ if and only if $\langle \rho_H \, , \, \mathcal{W} \rangle = 0$
 if and only if $\dim V^H = 0$.
\end{enumerate}

Furthermore, in the case $f_H \neq 0$, the left ideal
$\mathbb{Q}[G] \, f_H$ generated by the idempotent $f_H$ is a left
$\mathbb{Q}[G]-$module affording the representation $\mathcal{W}$
with multiplicity given by $\frac{\dim V^H}{m}$.
\end{thm}

\begin{proof}
The first statement is clear since $e_{\mathcal{W}}$ is central in
$ \mathbb{Q}[G] $ and both $p_H$ and $e_{\mathcal{W}}$ are
idempotents.

The second statement is obvious from the definition of $f_H$, and
it follows from it that $f_H$ is in $\mathcal{W}^H$.

Therefore $f_H = 0$ if $0 = \dim_{\mathbb{Q}} \mathcal{W}^H$.

We also have that
$$\dim_{\mathbb{Q}} \mathcal{W}^H = \langle \rho_H \, , \,
\mathcal{W} \rangle = \frac{\dim V^H}{m} \, \langle \mathcal{W} \,
, \, \mathcal{W} \rangle  = [L : \mathbb{Q}] \, \dim_{\mathbb{C}}
V^H
$$

\noindent where the second equality follows from Lemma
\ref{lem:genrho} and the third from (\ref{eq:WL}), thus proving
that $\langle \rho_H \, , \, \mathcal{W} \rangle = 0$ if and only
if $\dim_{\mathbb{C}} V^H = 0$.

We now show that if $\langle \rho_H \, , \, \mathcal{W} \rangle
\neq 0$ then $f_H \neq 0$ as follows.

 Since
  $$ f_H = \frac{1}{|H| |G|} \sum_{h \in H, \, g \in G}
  \chi_{\mathcal{W}}(g^{-1}) \ h \, g
  $$
 we compute the coefficient of $1$ of $f_H$ to be
  $$ \frac{1}{|H| |G|} \sum_{h \in H}
  \chi_{\mathcal{W}}(h) = \frac{1}{|G|} \langle \Res_H(\mathcal{W}) ,
  1_H\rangle_H = \frac{1}{|G|} \langle \mathcal{W}, \rho_H\rangle
  $$
and therefore $f_H \neq 0$ if $\langle \mathcal{W}, \rho_H\rangle
\neq 0$, which  finishes the proof of the third statement.

As for the last statement, the first result in Remark
\ref{rem:turull} applied to the case $F = \mathbb{Q}$ together
with the result of Lemma \ref{lem:genrho} imply that the left
ideal $\mathbb{Q}[G] \, f_H$ is equal to $\mathbb{Q}[G] \, p_H
\cap \mathbb{Q}[G] \, e_{\mathcal{W}}$ and it follows from Lemma
\ref{lem:genrho} that it is therefore equal to the direct sum of
$\frac{\dim_{\mathbb{C}} V^H}{m}$ minimal left ideals, each of
which affords the representation $\mathcal{W}$, thus finishing the
proof.
\end{proof}

\begin{cor}
\label{cor:turullgen} Continuing with the notation of Theorem
\ref{thm:f}, suppose that
 there exists a subgroup $H$ of $G$
such that
 $$ \langle \rho_H \, , \, V  \rangle_G = m_{\mathbb{Q}} (V) \ ,
 $$

\noindent or, equivalently, that there exist a subgroup $H$ of $G$
and a rational representation $\widetilde{\mathcal{W}}$ of $G$
with

 $$ \rho_H = \mathcal{W} \oplus \widetilde{\mathcal{W}} \ , \text{ and }
 \langle \mathcal{W} \, , \, \widetilde{\mathcal{W}} \rangle = 0 \ .
 $$

\medskip

Then $f_H$ is a primitive idempotent in $\mathbb{Q}[G] \,
e_{\mathcal{W}}$ commuting with every element $h$ in $H$ and
$\mathbb{Q}[G] \, f_H$ affords the representation $\mathcal{W}$.
\end{cor}

\section{An application: Group actions on abelian varieties}
\label{sec:abvar}

For a group $G$ acting on an abelian variety $A$, Lange and
Recillas \cite{L-R} have proved the following result, called the
\textit{isogeny decomposition} of $A$ with respect to $G$.

{\bf Theorem:} {\it Let $G$ be a finite group acting on an abelian
variety $A$. Let ${\mathcal{W}}_1, \ldots , {\mathcal{W}}_r$
denote the irreducible $\mathbb{Q}$-representations of $G$ and
$n_j:= dim_{D_j}({\mathcal{W}}_j)$ with $D_j :=
End_G({\mathcal{W}}_j)$ for $j = 1, \ldots , r$. Then there are
abelian subvarieties $B_1, \ldots , B_r$ such that each
$B_j^{n_j}$ is $G-$stable and associated to the representation
${\mathcal{W}}_j$, and also an isogeny}
\begin{equation}
\label{eq:isog}
 A \sim B^{n_1}_1 \times \ldots \times B^{n_r}_r \ .
\end{equation}

Observe that the integers $n_j$ in the theorem satisfy $n_j =
\frac{\dim V_j}{m_j}$, where $V_j$ is a complex irreducible
representation associated to ${\mathcal{W}}_j$ and $m_j =
m_{\mathbb{Q}} (V_j)$ is the Schur index of $V_j$.

To explain what the varieties $B_j$ are, observe that the action
of $G$ on the abelian variety $A$ induces an algebra homomorphism
$$
\rho : \mathbb{Q}[G] \to \End_{\mathbb{Q}}(A).
$$

If $e$ denotes any idempotent of the algebra $\mathbb{Q} [G]$,
define
$$
\Image (e) := \Image (\rho ( m e)) \subseteq A
$$
where $m$ is some positive integer such that $m e \in \mathbb{Z}
[G]$. $\Image (e)$ is an abelian subvariety of $A$, which
certainly does not depend on the chosen integer $m$.

The $B_j$ in the above theorem are then defined by $B_j = \Image
(p_1^j)$, where
$$ e_{{\mathcal{W}}_j} = p^j_1 + p^j_2 + \ldots + p_{n_j}^j
$$

\medskip

\noindent is any decomposition of the central idempotent
$e_{{\mathcal{W}_j}}$ as a sum of primitive orthogonal
idempotents.

Since there are, in general, many decompositions of
$e_{{\mathcal{W}}_j}$ of the form just given, it is clear that the
$B_j$ are far from unique, and defined only up to isogeny.

With the explicit construction of primitive orthogonal rational
idempotents $f$ given in Section \ref{sect:minimal} through
Theorem \ref{thm:kjs}, Remark \ref{rem:const} and Corollaries
\ref{cor:ks} and \ref{cor:fjs}, the $B_j$'s may now be obtained as
images of explicit idempotents, as in the following result.

\begin{prop}
Let $G$ be a finite group acting on an abelian variety $A$.

Then each $B_j$ in the isotypical decomposition (\ref{eq:isog}) of
$A$ with respect to $G$ may be chosen as the image of any of the
corresponding primitive rational idempotents $f$ constructed in
Corollary \ref{cor:fjs}.
\end{prop}

\subsection{The case of Jacobians}
\label{sec:Jac}

For the case of a $G-$action on the Jacobian variety $JW$ of a
curve $W$, we can be even more explicit.

The main new observation here is that every intermediate geometric
object specified by the cover $W \to W_G = W/G$, such as $JW_H =
J(W/H)$ or $P(W_H/W_N)$ for subgroups $H \subseteq N \subseteq G$,
has its own isogeny decomposition determined by an appropriate
combination of rational representation of $G$ of the form
$\rho_H$.

In fact, we will give the isotypical decomposition of any
intermediate Jacobian variety $JW_H$ and of any intermediate Prym
variety $P(W_H/W_N)$; we will also provide the rational
idempotents whose image is the corresponding $B_i^{s_i}$ inside
the given variety.

From now on we assume that the rational irreducible $\{
\mathcal{W}_1 , \ldots , \mathcal{W}r \}$ representations of $G$
are numbered so that $\mathcal{W}_1$ is the trivial
representation; therefore, the first factor in the isotypical
decomposition (\ref{eq:isog}) corresponding to a Jacobian $JW$
with $G-$action will always be $B_1^1 = JW_G$.

\begin{prop}
\label{prop:jac} Given a Galois cover $W \to W_G$, consider the
associated isotypical decomposition (\ref{eq:isog}) of $JW$ given
as
$$
JW \sim JW_G \times B^{\frac{\dim V_2}{m_2}}_2 \times \ldots
\times B^{\frac{\dim V_r}{m_r}}_r
$$

\noindent where $V_j$ is a complex irreducible representation
associated to $\mathcal{W}_j$ and $m_j = m_{\mathbb{Q}}(V_j)$ is
the Schur index of $V_j$.

Let $H$ be a subgroup of $G$ and denote by $\pi_H : W \to W_H$ the
corresponding quotient map.

Then the corresponding isotypical decomposition of $JW_H$ is given
as follows.

\begin{equation}
\label{eq:jacdec}
 JW_H \sim JW_G \times B^{\frac{\dim V_2^H}{m_2}}_2 \times \ldots
\times B^{\frac{\dim V_r^H}{m_r}}_r \ ,
\end{equation}

\medskip

\noindent where $V_j^H$ is the subspace of $V_j$ fixed by $H$.

Furthermore, considering images in $JW$ of rational idempotents as
defined at the beginning of Section \ref{sec:abvar}, and setting
$p_H = \frac{1}{|H|} \sum _{h \in H} h$ and $f_H^j = p_H \,
e_{\mathcal{W}_j}$ as in Theorem \ref{thm:f}, the following
equalities hold.

\begin{equation}
\Image (p_H) = \pi_H^{*} (JW_H)
\end{equation}

\noindent and if $\dim V_j^H \neq 0$ then
\begin{equation}
\Image (f_H^j) = B^{\frac{\dim V_j^H}{m_j}}_j \ .
\end{equation}
\end{prop}

\medskip

\begin{proof}
First, it is clear that $\Image (p_H)$ is the connected component
containing the origin of the subvariety $JW^H$ of $JW$ fixed by
$H$, and therefore equal to $\pi_H^{*} (JW_H)$, which is of course
isogenous to $JW_H$.

Secondly, it follows from Theorem \ref{thm:f} that the image of
each non-zero $f_H^j$ is contained in $\Image (p_H)$ and is equal
to $B^{\frac{\dim V_j^H}{m_j}}_j$.

Finally, the equality
$$ p_H = \sum_{j \in \{ 1 , \ldots , r\}} f_H^j
$$

\noindent induces the sought isogeny
$$ JW_G \times B^{\frac{\dim V_2^H}{m_2}}_2 \times \ldots
\times B^{\frac{\dim V_r^H}{m_r}}_r \to JW_H
$$

\noindent giving the isotypical decomposition of $JW_H$.
\end{proof}

\begin{rem}
The existence of non negative integers $s_j$ such that there is a
decomposition of the type given in the following corollary is
shown in \cite{L-R}; in the case when the group $G$ has only
absolutely irreducible representations, the $s_j$'s are shown to
be equal to $\dim V_j^H - \dim V_j^N$ in \cite{lr2}, using the
existence of an essentially unique $G-$invariant inner product on
each $\mathcal{W}_j = V_j$, which does not hold when the Schur
index is larger than one.
\end{rem}

We now give the isotypical decomposition of any intermediate Prym
variety.

\begin{cor}
\label{cor:genprym}

Given a Galois cover $W \to W_G$, consider the associated
isotypical decomposition of $JW$ given as
$$ JW \sim JW_G \times B^{\frac{\dim V_2}{m_2}}_2 \times \ldots
\times B^{\frac{\dim V_r}{m_r}}_r \ .
$$

Then for any subgroups $H \subseteq N \subseteq G$ the
corresponding decomposition of $P(W_H/W_N)$ is given as follows.

$$ P(W_H/W_N) \sim B^{s_2}_2 \times \ldots \times B^{s_r}_r
$$

\noindent where

$$ s_j =  \frac{\dim V_j^H}{m_j} - \frac{\dim V_j^N}{m_j} \ .
$$
\end{cor}

\begin{proof}
The previous Proposition gives the isotypical decompositions for
$JW_H$ and $JW_N$, but we also have the natural isogeny $JW_H \sim
JW_N \times P(W_H/W_N)$.

If we lift all these isogenies to the corresponding equalities
inside the tangent space to $JW$ at the origin, we obtain
equalities between sums of subspaces, which may be chosen to be
orthogonal among themselves with respect to the Hermitian product
induced by the natural polarization on $JW$ (see \cite{crr} for
similar arguments).

In other words, letting $\pi_H : W \to W_H$ denote the
corresponding cover map, the following decompositions are chosen
to be orthogonal decompositions.

\begin{equation}
\label{eq:pih} (d\pi_H^{*})_0 (T_0 JW_H) = (d\pi_G^{*})_0 (T_0
JW_G) \oplus T_0(B_2^\frac{\dim V_2^H}{m_2}) \oplus \ldots \oplus
T_0(B_r^\frac{\dim V_r^H}{m_r})
\end{equation}

\begin{equation}
\label{eq:pin} (d\pi_N^{*})_0 (T_0 JW_N ) = (d\pi_G^{*})_0 (T_0
JW_G) \oplus T_0(B_2^\frac{\dim V_2^N}{m_2}) \oplus \ldots \oplus
T_0(B_r^\frac{\dim V_r^N}{m_r})
\end{equation}

and

\begin{equation}
\label{eq:hn} (d\pi_H^{*})_0 (T_0 JW_H) = (d\pi_N^{*})_0 (T_0
JW_N) \oplus T_0(P(W_H/W_N))
\end{equation}

Now we can replace (\ref{eq:pin}) in (\ref{eq:hn}) and comparing
the result with (\ref{eq:pih}), the proof is finished.
\end{proof}

\begin{rem}
Note that $B_{\mathcal{W}}$ (for $\mathcal{W}$ different from the
trivial representation)  is always contained in some
$\pi_H^{*}(P(W_H/W_N))$, possibly in more that one.

This is because $B_{\mathcal{W}}$ is contained in
$\pi_H^{*}(JW_H)$ for any $H$ subgroup of $G$ such that $\langle
\rho_H , \mathcal{W} \rangle \neq 0$; in particular, this
condition is always satisfied for the trivial subgroup of $G$.
Moreover, for any $H$ satisfying $\langle \rho_H , \mathcal{W}
\rangle \neq 0$, there exist (at least) one subgroup $N$ of $G$
containing $H$ and such that $\langle \rho_N , \mathcal{W} \rangle
= 0$: for instance $N = G$, and for any such subgroup $N$,
$B_{\mathcal{W}}$ will be contained in $\pi_H^{*}(P(W_H/W_N))$.

In fact, in this case we know from the previous Corollary that it
will be contained with multiplicity $\frac{\dim
V^H}{m_{\mathbb{Q}}(V)}$, and from Proposition \ref{prop:jac} that
$B_{\mathcal{W}}^{\frac{\dim V^H}{m_{\mathbb{Q}}(V)}} \subseteq
\pi_H^{*}(P(W_H/W_N))$ may be obtained as the image of the
idempotent $f_H = p_H \,e_{\mathcal{W}}$.
\end{rem}

\begin{cor}
\label{cor:uniqueprym}

Given a Galois cover $W \to W_G$, consider the associated
isotypical decomposition of $JW$ given as
$$ JW \sim JW_G \times B^{\frac{\dim V_2}{m_2}}_2 \times \ldots
\times B^{\frac{\dim V_r}{m_r}}_r \ .
$$

Assume there exist subgroups $H \subset N \subseteq G$ and a
rational irreducible representation $\mathcal{W}$ of $G$ such that

\begin{equation}
\label{eq:uniqueprym}  \rho_H  = \mathcal{W} \oplus \rho_N \ .
\end{equation}

Then
$$ P(W_H/W_N) \sim B_{\mathcal{W}} \ .
$$

Conversely, if for some $\mathcal{W}$ rational irreducible
representation of $G$ there are subgroups $H \subset N \subseteq
G$ such that $P(W_H/W_N) \sim B_{\mathcal{W}}$, then
(\ref{eq:uniqueprym}) holds.
\end{cor}

\begin{proof}
Recall that the rational irreducible representations of $G$ are
denoted by $\{ \mathcal{W}_1 = 1, \ldots , \mathcal{W}_r \}$ and
$V_j$ denotes a complex irreducible representation associated to
$\mathcal{W}_j$.

For subgroups $H \subset N \subseteq G$, consider the non-negative
integers given by $s_j = \frac{\dim V_j^H}{m_j} - \frac{\dim
V_j^N}{m_j}$ for $j \in \{ 2 , \ldots , r\}$ and observe that it
follows from Lemma \ref{lem:genrho} that
$$ \rho_H - \rho_N = \bigoplus_{j=2}^r s_j \mathcal{W}_j \ .
$$

But it follows from Corollary \ref{cor:genprym} that $P(W_H/W_N)
\sim B_{\mathcal{W}}$ if and only if the unique $s_j$
corresponding to $\mathcal{W}$ equals one and all the other $s_j$
equal zero, thus finishing the proof.
\end{proof}

\begin{cor}
\label{cor:intersect}

Given a Galois cover $W \to W_G$, consider the associated
isotypical decomposition of $JW$ given as
$$ JW \sim JW_G \times B^{\frac{\dim V_2}{m_2}}_2 \times \ldots
\times B^{\frac{\dim V_r}{m_r}}_r \ .
$$

Assume there are subgroups $H$, $N_1$ and $N_2$ of $G$ such that
$H$ is a subgroup of $N_k$ for each $k$, and
 $$\rho_H =  \rho_{N_k} \oplus \mathcal{W} \oplus  \mathcal{W}_k
 $$

\noindent with $\mathcal{W}$, $\mathcal{W}_1$ and $\mathcal{W}_2$
rational representations of $G$ such that
$$ \langle \mathcal{W} , \mathcal{W}_j \rangle =
   \langle \mathcal{W}_1 , \mathcal{W}_2 \rangle = 0
$$

\noindent for $j = 1, 2$ and $\mathcal{W}$ irreducible.

Then, if $(X)^0$ denotes the connected component of $X$ containing
the origin, we have
\begin{equation}
B_{\mathcal{W}} \sim (P(W_H/W_{N_1}) \cap P(W_H/W_{N_2}))^0 \ .
\end{equation}
\end{cor}

It is clear that a similar result holds if there is a finite
number of subgroups $N_k$ containing $H$  such that the
representations $\rho_H - \rho_{N_k}$ all have a common rational
representation.

\begin{rem}
\label{rem:newva}

The previous corollaries can be used to give explicit geometric
decompositions of Jacobians with group actions in many cases, in
the manner illustrated with the examples given in the appendices.

Besides obtaining these isotypical decompositions, there also some
interesting new types of abelian subvarieties that emerge from the
picture as $B_{\mathcal{W}}$'s: the first examples worked out,
such as $D_p$ in \cite{ries}, $S_3$ in \cite{rr1}, the Klein
group, $D_4$, $A_4$ and $S_4$ in \cite{rr2} and $A_5$ in
\cite{sa1}, were all Prym or Jacobian varieties of intermediate
covers. In the case of the groups $D_n$ in \cite{crr} and also for
$S_5$ in \cite{sa2} and \cite{lr2}, a new type appeared: some of
the $B_{\mathcal{W}}$'s were nor Jacobian nor Prym varieties, but
(connected components of the origin of) intersections of Prym
varieties, of the form
$$ B_{\mathcal{W}} = \left( \bigcap_{k \in \{ 1, \ldots , s \}} P(W_H/W_{N_k}) \right)^{\circ}
$$

\noindent where the subgroups $H \subset N_k \subseteq G$ are as
in Corollary \ref{cor:intersect}. This last situation may  be
equivalently described (see \cite{lr2}) in the following way in
the case when $r = 2$. Consider the following diagram of curves
and covers

$$ 
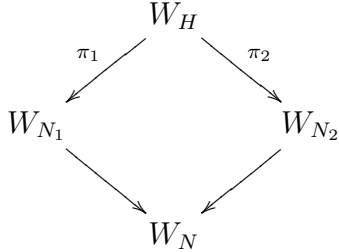
\xymatrix{
                & W_H \ar[dr]^{\pi_2} \ar[dl]_{\pi_1}            \\
 W_{N_1} \ar[dr]_{}  &  &     W_{N_2} \ar[dl]^{}\\
                & W_N        }
$$

\noindent with $N = \langle N_1 ,  N_2 \rangle $ and where $H$ ,
$N_1$ and $N_2$ satisfy the hypothesis of Corollary
\ref{cor:intersect}.

Then $B_{\mathcal{W}}$ is the (common) orthogonal complement of
$\pi_k^{*}(P(W_{N_k}/W_N))$ inside  $P(W_{H}/W_{N_l})$ for $k \neq
l$, where orthogonality is taken with respect to the Hermitian
product induced by the canonical Hermitian product induced by the
natural polarization on $JW_H$.

Now a further new type appears, by taking into account the Schur
index (all examples mentioned above have Schur index equal to one
for every complex irreducible representation). As may be seen in
the example in Appendix \ref{sect:otherexample}, there is  a
$B_{\mathcal{W}}$ that is not of any of the above types. However,
it is the image of an explicit rational idempotent and, as such,
it may be seen to be the orthogonal complement of $2$ Prym
varieties inside another Prym variety.

\end{rem}

\begin{rem}
\label{rem:newiso}

In the study of the decompositions of Jacobians in the case of the
special groups mentioned in the previous Remark, some interesting
isogenies have arisen. They were found with ad-hoc methods or with
restrictions on the acting representations, such as the group
having only absolutely irreducible representations.

It is clear now that these isogenies are due to equalities between
sums and differences of representations of $G$ of the type
$\rho_H$ for adequate subgroups $H$ of $G$; for instance,  the
equality
$$ \rho_S - \rho_R = \rho_X - \rho_G
$$

\noindent for $G =S_4$ and respective subgroups $X = S_3$, $R =
D_4$ and $S$ a Klein non normal subgroup, gives rise to an isogeny
between $P(W_S/W_R)$ and $P(W_X/W_G)$, which when studied in
detail (see \cite{rr2}) gives an alternate proof of the Recillas
trigonal construction \cite{rec}.

Since we now know that the representations $\rho_H$ for the
subgroups $H$ of $G$ determine (with no restrictions on $G$) the
isotypical decompositions of all the geometric intermediate pieces
in the full diagram of curves and covers of a curve with
$G-$action, it is likely that new interesting isogenies will
appear.
\end{rem}

\medskip

\appendix
\section{An example} \label{sect:example}

Consider the group of order $80$ given as follows

$$ G = \langle x, y : x^{20} \, , \,  y^8 \, , \, x^{10}y^4 \, , \, y^{-1}xyx^{-3}  \rangle
$$

\noindent and its irreducible representation $V$ of degree four
whose character $\chi_V$ is given in the following table on each
of the fourteen conjugacy classes of $G$.

\begin{center}
 \tiny{\begin{tabular}{|c|c|c|c|c|c|c|c|c|c|c|c|c|c|}
  \hline
                    &     &          &     &              &       &       &
           &       &              &         &       &        &  \\
           identity & $x$ & $x^{19}$ & $y$ & $x^{10} y^3$ & $x^2$ & $x y$ & $x^{11} y^3$
           & $y^2$ & $x^{10} y^2$ & $x y^2$ & $x^4$ & $x^{10}$ & $x^5$  \\
                    &     &          &     &              &       &       &
           &       &              &         &       &          &
           \\ \hline
                   &     &          &     &              &       &       &
           &       &              &         &       &          &\\
           $4$     & $k$ & $-k$     & $0$ & $0$          & $1$   & $0$   & $0$
           & $0$   & $0$          & $0$     & $-1$  & $-4$     & $0$ \\
                    &     &          &     &              &       &       &
           &       &              &         &       &          &  \\
  \hline
\end{tabular}
}
\end{center}

\bigskip

\noindent where $k = \sqrt{-5}$. Then
$$ K = \mathbb{Q}(k) \, \text{ and }\Gal (K/\mathbb{Q}) = \langle
\, \varphi(k) = -k \, \rangle \ .
$$

Furthermore, $V$ is defined over
$$ L = K(l) = \mathbb{Q}(\sqrt{10} + \sqrt{-2}) \, \text{ and }
\Gal(L/K) = \langle \, \tau(l) = -l \, \rangle \ ,
$$

\noindent where $l = \sqrt{-2}$, as seen from the following matrix
representation of $V$.

$$ x = \left(%
\begin{array}{cccc}
  \frac{2}{5} \sqrt{-5} & -\frac{1}{5} \sqrt{-5} & \frac{1}{5} \sqrt{-5} & -\frac{2}{5} \sqrt{-5} \\
  \frac{2}{5} \sqrt{-5} & \frac{1}{5} \sqrt{-5} & 0 & -\frac{1}{5} \sqrt{-5} \\
  0 & \frac{1}{5} \sqrt{-5} & \frac{2}{5} \sqrt{-5} & -\frac{2}{5} \sqrt{-5} \\
  \frac{1}{5} \sqrt{-5} & -\frac{1}{5} \sqrt{-5} & \frac{2}{5} \sqrt{-5} & 0 \\
\end{array}%
\right)
$$
 and
$$
 y = \left(%
\begin{array}{cccc}
  -1 + \sqrt{-2} & 0 & -1 - \sqrt{-2} & 1 \\
  -2 & 1 & -\sqrt{-2} &  \sqrt{-2} \\
  -1 &  \sqrt{-2} &  -\sqrt{-2} & -1 \\
  -1 & -1 &  1 - \sqrt{-2} & 0 \\
\end{array}%
\right)
$$

Then the Schur index $m_{\mathbb{Q}} (V) = 2$ and we should find
$\frac{\dim V}{m} = 2$ primitive orthogonal idempotents $k_1$ and
$k_2$ in $K[G] \, e_V$ and also two primitive orthogonal
idempotents $f_1$ and $f_2$ in $\mathbb{Q}[G] \, e_{\mathcal{W}}$
satisfying the following equalities.

\begin{multline} e_V  = k_1 + k_2 =
\frac{1}{20}(\text{identity} - x^{10})  \left( 4 \text{identity}
+ (-k) x \right. \\
 +  x^2 + (-k) x^3   \left. + (-1) x^4  +  x^6 + (-k) x^7
+ (-1) x^8 + (-k) x^9 \right)
\end{multline}

\begin{multline}
 e_{\mathcal{W}}  = e_V + e_{V^{\varphi}}  = f_1 + f_2 \\
 = \frac{1}{10}(\text{identity } - x^{10}) (4 \text{identity } + x^2 - x^4 + x^6 -
 x^{8})\ .
\end{multline}

We first compute the four primitive orthogonal idempotents
$\ell_j$ whose sum is $e_V$ as given by  (\ref{eq:special}). It
can be verified that $\ell_j$ and $\tau(\ell_j)$ are non
orthogonal.  We illustrate by writing down $\ell_1$, given as
follows.

\begin{multline}
\notag  \ell_1 = \frac{1}{100} (\text{identity} - x^{10})   \left(
5 (\text{identity} + y + x^2 + x^{14} y^3 + x^2 y^3 + x^4 y)
\right. \\
 + 10 (x^{12} y^2 + x^6 y^2) + (4 k) (x^{13} y^3 + x^{17} y)
+ (2 k) (x^{19} + x^{13})  + k (x^{15}+ x^{17}) \\
+ (k + k l)(x^{11} y^3 + x^{11} y) + (k -  kl) (x^3 y + x^{15} y^3) + (k - 2 k l) (x^{19} y + x^7 y^3)   \\
+ ( k + 2 kl) (x^{15} y + x^{19} y^3) + (3k + 3k l) (x^3 y^2)  + (3k - 3k l) (x^5 y^2) \\
+ (4 k + 2 k l) (x^{17} y^2) + (4 k - 2 k l) (x^{11} y^2)
+ (2 k l) (x^9 y^2) + (10 l) (x^{14} y^2)  \\
\left. + (5 + 5 l) (y^2  + x^{16} y^3 + x^6 y)  + (5 - 5 l) (y^3 +
x^{18} y^2 + x^{18} y) \right)
\end{multline}

Applying the construction detailed in Remark \ref{rem:const} we
obtain the new primitive orthogonal idempotents $u_1^1$ in $L[G]
\, \ell_1$ and $u_2^1$ in $L[G] \, \ell_2$ respectively, as
follows.

\begin{multline}
\notag u_1^1 = \frac{1}{80}(\text{identity} - x^{10}) \left(
4(\text{identity} + x^2) \right. + 3 (x^{18} + x^4)+ 6 (x^2 y +
x^{18} y^3) \\
+ 2 k (x^{15} + x^{13} y^3 + x^{17} y + x^{17}) + k (x^{13}+
x^{19}) + 2 l (x^{14} y^2) + 2 k l (x^{9} y^2) \\
+ (2+2l) (x^{10}y^2) + (2-2 l) (x^8 y^2) + (4+3 l) (x^4 y) + (4-3 l) (y)\\
+ (4+ l) (x^{16} y^3) + (4- l) y^3 + (1+4 l) (x^6 y)+ (1-4 l) (x^{18} y)\\
+ (1+3 l) (x^{14} y^3 + x^{16} y^2) + (1-3 l) (x^2 y^2 + x^2 y^3)
+ (2 k- k l) (x^{15} y + x^{15} y^3) \\
+ (2 k+ k l) (x^{11} y^3 + x^{19} y) + (k-2 k l) (x^{13} y)  +
(k+2 k l) (x y) \\
 \left. + (k+k l) (x^9 y^3 + x^7 y^2)  +
(k-k l) (x y^2 + x^{17} y^3) \right)
\end{multline}

\begin{multline}
\notag u_2^1 = \frac{1}{160}(\text{identity} - x^{10})
  \left(8(\text{identity}) + 4 k (x^5) \right. + l (x^{14} y^2) + (3 k l) (x^{17})\\
+ k l (x^{19} y^2 + x^{13} y^2) + (10-5 l) (x^{14}) + (8-2 l) (x^{10} y) + (8-6 l) (x^{10} y^3)\\
+ (8+5 l) (x^{14} y) + (8+3 l) (x^6 y^3) + (4+5 l) (x^6 + x^{12})
+ (4+4 l) (y^2) + (4+ l) (x^{18} y^2) \\
+ (4 k- k l) (x^3 y^3)  + (4 k+ k l) (x y^3 + x^7 y + x^{11})
+ (4 k+3 k l) (x^9 y) \\
+ (4 k-2 k l) (x^5 y^3 + x^5 y) + (12+ l) (x^{12} y) + (12- l)
(x^8 y^3) + (2+5 l) (x^4 y^3) \\
+ (2-5 l) (x^{8}) + (2+ l) (x^6 y^2) + (2- l) (x^{12} y^2) + (2-7
l) (x^{12} y^3 + x^8 y) \\
+ (2+9 l) (x^{16} y)
+ (2 k+ k l) (x^{11} y^2+ x^{17} y^2) + (2 k-k l) (x^{19}) \\
+ (2 k-3 k l) (x^7 y^3 + x^3 y) \left. + (2 k+3 k l) (x^{13} +
x^{11} y + x^{19} y^3) \right)
\end{multline}

\medskip

Then the $k_s$ are obtained applying Corollary \ref{cor:ks}; we
illustrate with $k_1$ as follows.
\begin{multline}
\notag k_1 = u^1_1 + \tau(u^1_1)
 = \frac{1}{40}(\text{identity} - x^{10})
\left( 4(\text{identity} + x^2 + x^4 y   \right.\\
+ y + x^{16} y^3 + y^3) + 3(x^{18} + x^4 + x^2 y + x^{18} y^3)\\
+ 2 (x^{10} y^2 + x^8 y^2) +  (x^6 y + x^{18} y + x^{14} y^3 + x^{16} y^2 + x^2 y^2 + x^2 y^3) \\
+ 2 k (x^{15} + x^{13} y^3 + x^{17} y + x^{17} + x^{15} y + x^{15} y^3 + x^{11} y^3 + x^{19} y)\\
 \left. + k (x^{13}+ x^{19} + x^{13} y + x y   + x^9 y^3 + x^7 y^2  + x y^2 + x^{17} y^3) \right)
\end{multline}

\medskip

Now the $f_s$ are obtained applying Corollary \ref{cor:fjs}, as
follows.
\begin{multline}
\notag f_1 = k_1 + \varphi(k_1) = u^1_1 + \tau(u^1_1) +
\varphi(u^1_1) + \varphi(\tau(u^1_1)) \\
= \frac{1}{20}(\text{identity} - x^{10})
\left( 4(\text{identity} + x^2 + x^4 y + y + x^{16} y^3 + y^3) \right.\\
+ x^6 y + x^{18} y + x^{14} y^3 + x^{16} y^2 + x^2 y^2 + x^2 y^3 \\
 \left. + 3 (x^{18} + x^4 + x^2 y + x^{18} y^3)+ 2 (x^{10} y^2 + x^8 y^2)   \right)
\end{multline}

\begin{multline}
\notag f_2 = k_2 + \varphi(k_2) = u_2^1 + \tau(u_2^1) +
\varphi(u^2_1) + \varphi(\tau(u^2_1)) \\
= \frac{1}{20}(\text{identity} - x^{10})
 \left(4 (\text{identity} + x^{10} y + x^{10} y^3 + x^{14} y + x^6 y^3))  \right.\\
+ 5 x^{14} + 6 (x^{12} y + x^8 y^3) + 2 (x^6 + x^{12} + y^2 + x^{18} y^2)\\
+ \left.  x^4 y^3 + x^{8} + x^6 y^2 + x^{12} y^2   + x^{12} y^3 +
x^8 y + x^{16} y \right)
\end{multline}

Note that the factor $(\text{identity} - x^{10})$ appearing in
both $f_1$ and $f_2$ indicates that their images are contained in
$P(W/W_{\langle x^{10}\rangle})$, and indeed one can verify that
$$ \langle \rho_{\{ 1 \}} , V \rangle = 4 \ \text{ and } \ \langle \rho_{\langle x^{10} \rangle} , V
\rangle = 0 \ .
$$

But another computation shows that
$$ \langle \rho_{H = \langle x y^{2} \rangle} , V
\rangle = 2
$$

\noindent and it follows from Proposition \ref{prop:jac} that
$B_{\mathcal{W}} = \Image (p_H e_{\mathcal{W}})$ is contained in
$\pi_H^{*}(JX_H)$ with multiplicity one, where $p_H
e_{\mathcal{W}}$ is given as follows.

\begin{equation}
p_H e_{\mathcal{W}}  = \frac{1}{20} (\text{identity } - x^{10})
(\text{identity } + x y^2) \left( 4 \text{identity } + x^2 - x^4 +
x^6 - x^{8}  \right)
\end{equation}

\bigskip

In order to give the isotypical decomposition of  a Jacobian $JW$
with $G-$action we first give some more information about the
group $G$.

The character table for $G$ is given as follows

{\tiny
 \begin{center}
\begin{tabular}{c|c|c|c|c|c|c|c|c|c|c|c|c|c|c|}
 &  &  &  &  &  &  &  &  &  &  &  &  &  &  \\
    & identity & $ x$ & $x^{19}$ & $y$ & $x^{10} y^3$ & $x^2$ & $x y$ &
  $x^{11} y^3$ & $y^2$ & $x^{10} y^2$ & $x y^2$ & $x^4$ &  $x^{10}$ & $x^5$
  \\
    &  &  &  &  &  &  &  &  &  &  &  &  &  &  \\  \hline
    $\chi_{V_1}$ & 1 & 1 & 1 & 1 & 1 & 1 & 1 & 1 & 1 & 1 & 1 & 1 & 1 & 1 \\
     $\chi_{V_2}$ & 1 & -1 & -1 & -1 & -1 & 1 & 1 & 1 & 1 & 1 & -1 & 1 & 1 & -1
    \\
     $\chi_{V_3}$ & 1 & -1 & -1 & 1 & 1 & 1 & -1 & -1 & 1 & 1 & -1 & 1 & 1 & -1
    \\
     $\chi_{V_4}$ & 1 & 1 & 1 & -1 & -1 & 1 & -1 & -1 & 1 & 1 & 1 & 1 & 1 & 1 \\
     $\chi_{V_5}$ & 1 & -1 & -1 & $- w_4$ & $w_4$ & 1 & $w_4$ & $- w_4$ & -1 & -1 & 1 & 1 & 1 & -1
    \\
     $\chi_{V_6}$ & 1 & -1 & -1 & $w_4$ & $- w_4$ & 1 & $- w_4$ & $w_4$ & -1 & -1 & 1 & 1 & 1 & -1
    \\
     $\chi_{V_7}$ & 1 & 1 & 1 & $- w_4$ & $w_4$ & 1 & $- w_4$ & $w_4$ & -1 & -1 & -1 & 1 & 1 & 1
    \\
     $\chi_{V_8}$ & 1 & 1 & 1 & $w_4$ & $- w_4$ & 1 & $w_4$ & $- w_4$ & -1 & -1 & -1 & 1 & 1 &
    1\\
     $\chi_{V_9}$ & 2 & 0 & 0 & 0 & 0 & -2 & 0 & 0 & -2 $w_4$ & 2 $w_4$ & 0 & 2 & -2 & 0
    \\
     $\chi_{V_{10}}$ & 2 & 0 & 0 & 0 & 0 & -2 & 0 & 0 & 2 $w_4$ & -2 $w_4$ & 0 & 2 & -2 & 0 \\
     $\chi_{V_{11}}$ & 4 & 1 & 1 & 0 & 0 & -1 & 0 & 0 & 0 & 0 & 0 & -1 & 4 & -4 \\
     $\chi_{V_{12}}$ & 4 & -1 & -1 & 0 & 0 & -1 & 0 & 0 & 0 & 0 & 0 & -1 & 4 & 4 \\
     $\chi_{V_{13}}$ & 4 & $-k$ & $k$ & 0 & 0 & 1 & 0 & 0 & 0 & 0 & 0 & -1 & -4 & 0 \\
     $\chi_{V_{14}}$ & 4 & $k$ & $-k$ & 0 & 0 & 1 & 0 & 0 & 0 & 0 & 0 & -1 & -4 & 0 \\
  \hline
\end{tabular}
\end{center}
}

\bigskip

\noindent where $w_n$ denotes a primitive n-th root of unity and
$k = \sqrt{-5} = w_{20}+w_{20}^9-w_{20}^{13}-w_{20}^{17}$.

\medskip

Note that the irreducible rational representations of $G$ are
$V_j$, for $j$ in $\{ 1, 2 , 3, 4 , 11 , 12\}$, and $V_5 \oplus
V_6$, $V_7 \oplus V_8$, $V_9 \oplus V_{10}$, $2(V_{13} \oplus
V_{14})$. The representation $V$ given at the beginning of this
section is $V_{13}$.

\bigskip

Consider the following subgroups of $G$: $H_1 = \langle x^2 , x y
\rangle$, $H_2 = \langle x^2 , y \rangle$, $H_3 = \langle y^2 , x
\rangle$, $H_4 = \langle x^2 , x y^2 \rangle$, $H_5 = \langle x
\rangle$, $H_6 = \langle x^4 , x y^2 \rangle$, $H_7 = \langle x^3
y^2 , x y \rangle$, $H_8 = \langle x y\rangle$, $H_9 = \langle x
y^2 \rangle$, $H_{10} = \langle x y^2, x^{10} \rangle$ . Then we
have the following result.

\begin{thm}
Let $W$ denote a compact Riemann surface with $G-$action and set

$$ B_1 = P(W_{H_8}/W_{H_7}) \cap P(W_{H_8}/W_{H_1}) \ , \text{ and}
$$

$$ B_2 = P(W_{H_9}/W_{H_{10}}) \cap P(W_{H_9}/W_{H_6}) \ .
$$

\medskip

Then the isotypical decomposition of the Jacobian variety $JW$
with respect to $G$ is given as follows
\begin{multline}
\notag JW \backsim JW_G \times P(W_{H_1}/W_G) \times
P(W_{H_2}/W_G)
\times P(W_{H_3}/W_G) \\
\times P(W_{H_4}/W_{H_3}) \times P(W_{H_5}/W_{H_3})  \times
 P(W_{H_6}/W_{H_4})^2 \\
\times  P(W_{H_7}/W_G)^4 \times  B_1^4 \times  B_2^2
\end{multline}

\noindent with associated rational representations on the right
hand side given as follows.
\begin{multline}
\notag {} \ \ V_1 \oplus V_2 \oplus V_3 \oplus V_4 \\
\oplus (V_5 \oplus V_6) \oplus (V_7 \oplus V_8) \oplus
(V_9 \oplus V_{10}) \\
\oplus V_{12} \oplus V_{11} \oplus 2(V_{13}\oplus V_{14})
\end{multline}
\end{thm}

\begin{proof}
The proof is obtained by first computing the corresponding
representations $\rho_H$ for all (conjugacy classes of) subgroups
 $H$ of $G$, then looking for each irreducible rational representation
$\mathcal{W}$ in terms of the form $\rho_H - \rho_N$, for $H
\subset N \subseteq G$, and finally applying either Corollary
\ref{cor:uniqueprym} or Corollary \ref{cor:intersect}.

The first two steps can be implemented very quickly with a simple
computer program (see \cite{ar} for a GAP version).
\end{proof}

\section{Another example: a quasi-Prym}
\label{sect:otherexample}

Consider the following group, a $\mathbb{Z}/3\mathbb{Z}$ extension
of the quaternion group of order $8$,
$$ G = \langle x^4 , y^4 , z^3, y^{-1} x y x, z^{-1} x z y^{-1}, z^{-1} y z (xy)^{-1}
\rangle\ .
$$

Then $G$ has a complex irreducible representation $V$ of degree
two such that $K = \mathbb{Q}$ and $L = \mathbb{Q}(w_3)$, with
$w_3$ a primitive cubic root of $1$. Therefore $m_{\mathbb{Q}}(V)
= 2$ and $L \otimes_{\mathbb{Q}} \mathcal{W} = 2 \, V$.

Since $\frac{\dim V}{m_{\mathbb{Q}} (V)} = 1$, the factor
corresponding to $\mathcal{W}$ in the isotypical decomposition of
a Jacobian $JW$ with $G-$action will be of the form
$B_{\mathcal{W}}^1$.

It is clear that in this case
 \begin{align*} e_{\mathcal{W}} & = e_{V}  \\
                 & = \frac{1}{12} (\id -x^2) \, (2\id -(z + xz + yz + xyz)-(z^2+x^3z^2+y^3z^2+x^3yz^2))
\end{align*}

\noindent is a primitive rational idempotent.

Furthermore, it is known (c.f., \cite[pp. 177]{berkovich}) that
when $m_{\mathbb{Q}} (V) = \dim V$ then
 $$ \langle \rho_H , V \rangle \neq 0 \text{ if and only if } H =
 \{ 1 \} \ .
 $$

\medskip

This fact together with the explicit expression for
$e_{\mathcal{W}}$ given above imply that
$$B_{\mathcal{W}} = \Image(e_{\mathcal{W}}) \subset P(W/W_{\langle
 x^2 \rangle}) \cap P(W/W_{\langle z \rangle}) \ .
$$

One can also verify the following identities
\begin{align}
\label{eq:rhox2}
\rho_{1} - \rho_{\langle  x^2 \rangle} & = \mathcal{W} \oplus 2\mathcal{W}_1\\
\notag \rho_{1} - \rho_{\langle z \rangle} & = \mathcal{W} \oplus
\mathcal{W}_1 \oplus \mathcal{W}_2 \\
\notag \rho_{\langle z \rangle} - \rho_{\langle z , x^2 \rangle} &
= \mathcal{W}_1
\end{align}

\noindent where $\mathcal{W}_1$ and $\mathcal{W}_2$ are rational
representations of $G$ that are mutually disjoint and also
disjoint from $\mathcal{W}$.

Furthermore $\mathcal{W}_1$ is irreducible, and it follows from
Corollary \ref{cor:genprym} that the isotypical decomposition of
the (connected component of the origin of the)  intersection of
the two Prym varieties is the following.

$$ \left( P(W/W_{\langle  x^2 \rangle}) \cap P(W/W_{\langle z
\rangle}) \right)^0 \sim B_{\mathcal{W}} \times B_{\mathcal{W}_1}
\ .
$$

Some further computations allow us to exclude other possibilities
and therefore show that $B_{\mathcal{W}}$ is not any of the known
types mentioned in Remark \ref{rem:newva}; that is, it cannot be
isolated neither as a Prym variety nor as the intersection of Prym
varieties nor as the orthogonal complement of a Prym inside
another Prym variety.

However, the relations (\ref{eq:rhox2}) show that
$B_{\mathcal{W}}$ is the orthogonal complement of
$B_{\mathcal{W}_1}^2$, which is isogenous to $P(W_{\langle z
\rangle}/W_{\langle z , x^2 \rangle})^2$, inside $P(W/W_{\langle
x^2 \rangle})$.

\bibliographystyle{amsplain}

\end{document}